\def\sqr#1#2{{\vcenter{\hrule height.#2pt
     \hbox{\vrule width.#2pt height#1pt\kern#1pt
     \vrule width.#2pt}
     \hrule height.#2pt}}}
\def\square{\mathchoice\sqr{5.5}4\sqr{5.0}4\sqr{4.8}3\sqr{4.8}3}
\def\qed{\hskip4pt plus1fill\ $\square$\par\medbreak}
\def\cD{{\cal D}}
\def\CC{{\rm\kern.24em\vrule width.02em height1.4ex depth-.05ex\kern-.26em C}}
\def\RR{\,\rm{\vrule width.02em height1.55ex depth-.07ex\kern-.3165em R}}

\def\wt{{ wt}\,}
\def\cx#1{{\CC}^{#1}}
\def\vf#1{{\cal A}^{(#1)}}
\magnification=\magstep1
\rightline{February 1, 1992}
\bigskip\bigskip
\centerline{\bf Domains in $\cx {n+1}$ with}
\centerline{\bf Noncompact Automorphism Group. II}
\medskip
\centerline{Eric Bedford and Sergey Pinchuk}
\bigskip\bigskip
\centerline{\bf\S1. Introduction}\bigskip
We consider here smoothly bounded domains with
noncompact automorphism groups.  Examples of such domains may be obtained as
follows.  To the variables $z_1,\dots,z_n$ we assign weights
$\delta_1,\dots,\delta_n$ with $\delta_j=1/2m_j$ for $m_j$ a positive integer. 
If $J=(j_1,\dots,j_n)$ and $K=(k_1,\dots,k_n)$ are multi-indices, we set $wt\,
(J)=j_1\delta_1+\dots+j_n\delta_n$ and $wt\,(z^J\bar z^K)=wt\, J+wt\, K$.  We
consider real polynomials of the form $$p(z,\bar z)=\sum_{\wt J=\wt K={1\over
2}}a_{JK}z^J\bar z^K.\eqno(1.1)$$ 
The reality condition is equivalent to
$a_{JK}=\bar a_{KJ}$.  The balance of the weights, i.e. $\wt J=\wt K$ for each
monomial in $p$, implies that the domain $$G=\{(w,z_1,\dots,z_n)\in{\cx{}}\times
\cx{n}: |w|^2+p(z,\bar z)<1\}\eqno(1.2)$$  is invariant under the $T^2$-action 
$$(\varphi,\theta)\mapsto
(e^{i\varphi}w,e^{i\delta_1\theta}z_1,\dots,e^{i\delta_n\theta}z_n)\eqno(1.3)$$
In addition, the weighted homogeneity of $p$ allows the  transform
$(w,z)\mapsto(w^*,z^*)$ defined by 
$$w=\left(1-{iw^*\over
4}\right)\left(1+{iw^*\over 4}\right)^{-1},\quad z_j= z_j^*\left(1+{iw^*\over
4}\right)^{-2\delta_j}  \eqno(1.4)$$  to map the domain $G$ biholomorphically
onto  
$$D=\{(w,z_1,\dots,z_n)\in\cx{}\times\cx n: \Im m\,w+p(z,\bar
z)<1\}.\eqno(1.5)$$ This unbounded representation of $G$ shows that the
automorphism group is noncompact, since it is invariant under translation in the
$\Re e\,w$-direction.  Since $p$ is homogeneous, there is also a 1-parameter
family dilations; so with (1.3), $Aut(G)$ has dimension at least 3.

Presumably, every smoothly bounded domain $\Omega\subset\cx {n+1}$ with
noncompact automorphism group is equivalent to a domain of the form (1.2).  
Several papers have been written on the general subject of describing a such
domains, starting with the work of R.~Greene and S.~Krantz (see [BP] for a more
complete list of references).  Our approach to this problem may be thought of as
involving two steps.  The first is to show that $\Omega$ is biholomorphically
equivalent to one of a more special class of model domains $D$.  Such a model
domain would have a nontrivial holomorphic vector field.  The second step is to
transport this vector field back to $\Omega$, analyze such a vector field
tangent to $\partial\Omega$ at the parabolic fixed point, and use this
information determine the original domain.

The first progress in this direction was obtained for the first step. 
S.~Frankel [F] (see also Kim [K]) showed that if $\Omega$ is a (possibly
nonsmooth) convex domain, then $\Omega$ has an unbounded representation which
is independent of the variable $u$.  This particular model domain,
however, is not well enough behaved for us to carry out step 2 of our program.
In [BP] we used the approach described above in the case where the Levi form of
$\partial\Omega$ had at most one zero eigenvalue, and we obtain model domains
with polynomial boundaries.  

Here we
consider the case of convex domains, which makes the scaling arguments of [BP]
easier.  Thus we are able to complete step 1 more quickly and make
the transition to step 2.  Now we focus on the second step and show that if a
domain has a vector field of positive weight (defined below), then the domain
must be of the form (1.2).  We note that our main result, Theorem 3.7, is proved
for domains considerably more general than convex.   

By
{\it holomorphic vector field} we mean a vector field of the form
$H=\sum_{j=0}^n H_j{\partial\over\partial z_j}$, where $H_j$ is holomorphic on
$\Omega$,  and we use the notation $z_0=w$.  If the coefficients of $H$
extend smoothly to $\bar\Omega$, and  if the real vector field $\Re e\,H={1\over
2}(H+\bar H)$ is tangent to $\partial\Omega$, then we will say that $H$ is {\it
tangent} to $\partial\Omega$.  If $H$ is a holomorphic vector field which is
tangent to $\partial\Omega$, then the exponential $\exp(tH)$ defines an
automorphism of $\Omega$ for any $t\in \RR$.  

We assign weight 1 to the
variable $w=z_0$, i.e. $\delta_0=1$.  Thus we set
$wt\,(z^J{\partial\over\partial z_k})=wt\,(J)-\delta_k$ for any multi-index
$J=(j_0,\dots,j_n)$, and $0\le k\le n$.  If $H$ is a holomorphic vector field
(not necessarily homogeneous), we let $wt\, H$ denote the smallest weight of a
nonzero homogeneous term in the Taylor expansion of $H$ at 0.  If $H$ vanishes
to infinite order, we set $wt\, H=\infty$.  The one-parameter subgroup of
translations in the $\Re e\,w$-direction (automorphisms of $D$), are generated by
the vector field $2{\partial\over\partial w}$, which has a fixed point at
$\infty$.  If $p\ge0$, then $w\mapsto -1/w$ is well-defined on $D$, and the
vector field $2{\partial\over\partial w}$ is taken into
$2w^2{\partial\over\partial w}$, which has a fixed point at 0 and weight 1
there.

Let us fix a point $(w^0,z^0)\in\partial\Omega$.  After a translation and
rotation of coordinates, we may assume that $(w^0,z^0)=(0,0)$ and that near
$(0,0)$ we have $\Omega=\{v+f(u,z,\bar z)\}$ with $f(0,0)=\nabla f(0,0)=0$. An
assignment of weights $\delta_1,\dots,\delta_n$ is {\it admissible} if there is a
homogeneous polynomial $p(z,\bar z)$ (consisting of all the terms in the Taylor
expansin of $f$ weight 1)  such that all the terms in the Taylor expansion of
$f-p$ have weight strictly greater than 1.  In this case we define
$$\Omega_{hom}=\{v+p(z,\bar z)<0\}$$ 
to be the {\it homogeneous model} of $\Omega$ at (0,0).  (In general the
homogeneous model depends on the choice of weights.)  

\proclaim Theorem 1.  Let $\Omega$ be a domain with smooth boundary, and
suppose that there is an assignment of weights at a point
$(0,0)\in\partial\Omega$ such that  in the homogeneous model $p\ge0$, and $p$
does not vanish on any complex variety passing through 0.  If there is a
tangential holomorphic vector field $H$ for $\Omega$ with $H(0)=0$ and $0<\wt
H<\infty$, then $\Omega$ is biholomorphically equivalent to a domain of the form
(1.2).

A tangential holomorphic vector field $H$ generates a 1-parameter group of
automorphisms via the exponential map $t\mapsto\exp(tH)$.  The hypothesis that
$wt\,H>0$ in Theorem 1 serves to eliminate the possibility that $H$ generates a
local translation at (0,0) (in which case $wt\,H<0$) or a local dilation or
rotation (in both cases $wt\,H=0$). Our motivation in proving Theorem 1 is that
with the addition of the results and techniques of [BP] we obtain the
following.

\proclaim Theorem 2.  Let $\Omega\subset\cx{n+1}$ be a bounded, convex set with
smooth, finite type boundary.  If $Aut(\Omega)$ is noncompact, then $\Omega$ is
biholomorphically equivalent to a convex domain of the form (1.2).

Theorem 1 does not make any assumption of pseudoconvexity.  Our
arguments need some sort of convexity hypothesis, however, to pass from the
noncompactness of $Aut(\Omega)$ to the existence of a vector field $H$.

It is worth noting that Theorem 2 gives a classification, up to
biholomorphism, of the convex, finite type domains with noncompact automorphism
groups.  For by (1.3), a domain of the form (1.2) is a Cartan domain, and by [KU]
it follows that two domains of the form (1.2) are biholomorphically equivalent if
and only if they are linearly equivalent.  

Most of the following paper is devoted to the analysis of homogeneous vector
fields tangent to $\partial\Omega$.  Only certain vector fields can arise.  In
\S2 we show that there can be no tangent vector field of the form (2.5), i.e.\
independent of the variable $z_0=w$.  The next major case to analyze is the
case of vector fields of positive weight.  This is done in \S3.  In \S4 we show
how to apply these algebraic results to domains with noncompact automorphism
groups.  For the most part, this is a reiteration of results of [BP].  In the
convex case, we have more flexibility in our normal family arguments.  Thus we
give Lemma 4.1, which makes the transition between steps one and two more
natural and understandable than the treatment in [BP, Lemma 7].

\bigbreak

\centerline{\bf \S2. Holomorphic Tangent Vector Fields}
\bigskip

A vector field on $\cx n$ may be written as $Q=\sum q_j(z){\partial
\over\partial z_j}$, and $Q$ is homogeneous of weight $\mu$ if $\wt
q_j=\mu-\delta_j$ holds for $1\le j\le n$.  An {\it integral
curve} of $Q$ is a holomorphic function $\varphi:\cD\to\cx n$ for some domain
$\cD\subset\cx{}$ such that $\dot\varphi(t)=Q(\varphi(t))$ for all $t\in\cD$. 
If $\cD$ is a maximal domain of analyticity of an integral curve $\varphi$,
then the image $\varphi(\cD)$ is an {\it complex orbit} of $Q$.  

A complex orbit is necessarily unbounded unless it is a
constant (i.e. a critical point).  For if $\varphi:D\to\cx n$ is an integral
curve, and if $|\varphi(t)|<M$ on $\cD$ then for any point $t_0\in \cD$, the
solution of $\dot\varphi(t)=Q(\varphi(t))$ may be analytically extended to a disk $\{|t-t_0|<r\}$ where
$r$ is independent of $t_0$.  Thus $\varphi$ extends to an integral curve
$\varphi:\cx{}\to\cx n$.  Since $\varphi$ is bounded, it must be constant,
and so $\varphi(\cD)$ must be a critical point of $Q$.

For $\tau\in\cx{}$ we
define the dilation $D_\tau:\cx n\to\cx n$ by
$$D_\tau(z)=(\tau^{\delta_1}z_1,\dots,\tau^{\delta_n}z_n),$$ where
$\tau^{\delta_j}$ denotes an arbitrary but fixed choice of fractional power.  If
$\varphi$ is a solution of the vector field $Q$, then so is
$\varphi_\tau(t)=D_\tau(\varphi(\tau^\mu t))$.  Thus if $S$ is an orbit of $Q$,
then so is $D_\tau(S)$.

We start with some observations about homogeneous vector fields on $\cx n$.

\proclaim Lemma 2.1.  Let $Q$ be a nonzero homogeneous vector field of weight
$\mu\ne0$, and set $A:=\{Q=0\}$.  If $A=\{0\}$, then there is an orbit $S$ of
$Q$ which contains 0 in its closure.  \vskip0pt\indent
More generally, let $\Gamma\subset\cx n$ be a subvariety such that
$A\cap\Gamma=\{0\}$ and $\Gamma$ is invariant under $D_\tau$ for all
$\tau\in\cx{}$.  If $Q$ is tangential to $\Gamma$ at all regular points, then
there exists an orbit $S$ of $Q$ with $S\subset \Gamma$ and which contains 0 in
its closure.

\noindent{\sl Proof.}  We assume first that $A=\{0\}$.  Let $S^\prime$ be a
nontrivial orbit of $Q$ passing through a point $b_0$, and take a sequence
$a_k\in S^\prime$, $|a_k|\to\infty$.  For each $k$, we choose a dilation
$D_\epsilon$ with $\epsilon=\epsilon_k$ chosen such that $|D_\epsilon(a_k)|=1$. 
For each $k$, $S_k:=D_\epsilon(S^\prime)$ is an orbit of $Q$ and contains the
point $b_k:=D_\epsilon(b_0)$.  Since $b_k\to 0$ as $k\to\infty$, we may assume
that $|b_k|<1$.  Thus there exists a point $a_k\in S_k$  with $|a_k|=1$. 
Passing to a subsequence, we may assume that $a_k\to a_0$ with $|a_0|=1$.  Since
$Q(a_0)\ne0$, there is a nontrivial orbit $S$ passing through $a_0$, and it is
a basic property of systems of ordinary differential equations that the orbits
$s_k$ converge to the orbit $S$.  Since the $b_k\in S_k$ converge to 0, it
follows that 0 is in the closure of $S$.

The proof of the second assertion is identical; we just observe that $S_k$ and
$S$ remain inside $\Gamma$.\qed

\proclaim Lemma 2.2.  Let $Q$ be a nonzero homogeneous vector field of
weight $\mu\ne0$ which is not identically zero.  Then there is a nontrivial set
$\Sigma$ which is invariant under $D_\tau$ for all $\tau\in\cx{}$ such that any
solution $f$ of $Qf=0$ satisfies $f=f(0)$ on $\Sigma$.

\noindent{\sl Proof.}   Let us fix a
point $c$ with $Q(c)\ne0$, and let $S$ denote the orbit passing through $c$. 
Let $\Sigma$ denote the set of all limit points of sequences $\{a_\nu\}$ with
$a_\nu\in D_{\tau_\nu}S$ and $\tau_\nu\to 0$.  We note that
$\Sigma\cap\{|z|=1\}\ne\emptyset$.  It follows that $\Sigma$ is invariant under
$D_\tau$.  Further if $a_0\in\Sigma$, then there exists a sequence $a_\nu\in
D_{\tau_\nu}S$ with $a_\nu\to a_0$ and $\tau_\nu\to0$.  Thus
$b_\nu:=D_{\tau_\nu}c\to0$.  If $Qf=0$, then $f$ constant on $D_{\tau_\nu}S$,
so $f|_{D_{\tau_\nu}S}=f(b_\nu)$.  Thus
$f(a)=\lim_{\nu\to\infty}f(b_\nu)=f(0)$.  Since this holds for all
$a\in\Sigma$, $f|_\Sigma=f(0)$.\qed 

Any polynomial of $z$ and $\bar z$ may be
written in the form $p=\sum_{A,B} c_{A,B}z^A\bar z^B $.  We define the {\it
signature} $\delta(A,B)=\wt A -\wt B$.  Thus we may write $p$ as a sum
of monomials of fixed signature
$$p=p^{(\nu_{-L})}+p^{(\nu_{-L+1})}+\dots+p^{(\nu_{L-1})}+p^{(\nu_L)},\eqno(2.1)
$$ where $\nu_{-l}=-\nu_l$, $\nu_0=0$, and $\nu_l<\nu_{l+1}$.  Thus
$p^{(\nu_{-l})}=\overline{p^{(\nu_l)}}$.  We may define holomorphic functions
$f_{\nu,B}$ by summing first over the indices $A$:
$$p^{(\nu)}=\sum_{\delta(A,B)=\nu}c_{A,B}z^A\bar z^B=\sum_B f_{\nu,B}(z) \bar
z^B.\eqno (2.2)$$ We say that $p$ is {\it balanced} if $p=p^0$.  A
homogeneous polynomial of weight 1 is balanced if and only if it has the form
(1.1).  We will assume throughout the rest of this paper that $p$ is
homogeneous of weight 1.  Thus each of the functions $f_{\nu, B}$ is homogeneous
of weight $(\nu+1)/2$.

\proclaim Lemma 2.3.  Let $Q$ be a homogeneous vector field with
$Q(0)=0$, and suppose that $S$ is a nontrivial orbit of $Q$ with $0\in\overline
S$, and  such that $Q\ne0$ on $S$.  If $p^{(\nu)}$ is in the form (2.2) and if 
$Q^\alpha p^{(\nu)}=0$ on $\cx n$, then $p^{(\nu)}|_{S}=0$.

\noindent{\sl Proof.}  Since the vector field $Q$ is tangent to $S$, it follows
that $(Q^\alpha p^{(\nu)})|_S=(Q^\alpha|_S)(p^{(\nu)}|_S)$.  Thus from the
equation
$$Q^\alpha p^{(\nu)}=\sum Q^\alpha f_{\nu,B}(z)\bar z^B=0$$
we deduce that $Q^\alpha f_{\nu,B}=0$ on $\cx n$ for all $\nu$ and $B$, and
thus this holds on $S$.  Since $Q\ne0$ on $S$,  $Q^{\alpha-1}f_{\nu,B}$ is
constant on $S$.  But since $Q(0)=0$, we have $Q^{\alpha-1}f_{\nu,B}(0)=0$, and
since $0\in\bar S$ it follows that $Q^{\alpha-1}f_{\nu,B}|_S=0$.  Proceeding in
this way, we have $Q^\alpha f_{\nu,B}|_S=\dots= Qf_{\nu,B}|_S=f_{\nu,B}|_S=0$.
Since this holds for all $\nu$ and $B$, we conclude that $p^{(\nu)}=0$ on $S$.
\qed

We
let $\vf{\mu}$ denote the set of holomorphic vector fields which are
homogeneous of weight $\mu$ and which are tangential to the domain
$\Omega_{hom}$.
With the notation
$$Q=q_0{\partial\over \partial w} +  \sum_{j=1}^n q_j{\partial\over \partial
z_j}$$
the tangency condition is given by  
$$\Re
e\left(-{i\over2}q_0+\sum_{j=1}^n q_j {\partial p\over\partial
z_j}\right)=0\eqno(2.3)$$ for all $(w,z)\in\partial\Omega_{hom}$.  

Without loss of generality we may assume that
$$p=\sum_{A,B} c_{A,B}z^A\bar z^B {\rm\ 
contains\ no\ holomorphic\ (or\ antiholomorphic)\ monomials},\eqno(2.4)$$
i.e. neither multi-index $A$ or $B$ in the summation is equal to $(0,\dots,0)$.
Since $p$ is real, this is equivalent to the condition $f_{\nu,B}(0)=0$ for all
$f_{\nu,B}$ in (2.2).

Our first step will be to eliminate tangent vector fields of the special form
(independent of the variable $w$) 
$$Q=\sum_{j=1}^n
q_j(z_1,\dots,z_n){\partial\over\partial z_j}\eqno(2.5)$$ from the context of
Theorems 1 and 2.  Such a vector field can belong to $\vf\mu$ in the degenerate
case $Qp=0$.  Vector fields of the form (2.5) arise, too, as rotations $\sum
c_{ij}z_i{\partial\over\partial z_j}$ in the case of weight zero.  Another
possibility with positive weight is as follows.

\medskip
{\sl Example.}  Let $f=z_1^3z_2^2$, $g=z_1z_2$, $p=2\Re e\,f\bar g$, and
$Q=iz_1^2z_2(2z_1{\partial\over\partial z_1}- 3z_2{\partial\over\partial z_2}
)$.  Then $Qp=-if\bar f$, and $\Re e\,Qp=0$.

In order to eliminate the possibility of nonzero vector fields of the form
(2.5) with weight $\mu\ne0$,  we will make the hypothesis 
$$\{p=0\}{\rm\ contains\ no\ nontrivial\ complex\
manifold}.\eqno(2.6)$$

We will say that $(\tilde z_1,\dots,\tilde z_n)$
is a {\it weighted change of coordinates}  
if $\tilde z_j$ is a homogeneous polynomial of $(z_1,\dots,z_n)$, and its
weight is equal to $wt\,z_j$. This leads to a new polynomial $\tilde p$,
defined by $\tilde p(\tilde
z)=p(z)$.  It is evident that $\tilde p$ satisfies  (2.4) and (2.6) if and only
if $p$ does. Similarly, $\tilde p$ is balanced if and only if $p$ is.  

\proclaim Lemma 2.4.  If $Q\in\vf\mu$, then $Q(0)=0$.  

\noindent{\sl Proof.}  For
otherwise, $q_j=c\ne0$ for some $j$.  In this case we may make a homogeneous
change of coordinates to bring $Q$ into the form ${\partial\over\partial z_j}$. 
This means that $p$ is independent of the variable $x_j$.  However, setting
$z_k=0$ for $k\ne j$, we must have a nontrivial homogeneous polynomial of weight
$\mu$ in the variable $z_j$ alone, by (2.6).  On the other hand, the only
polynomial independent of $x_j$ is $(i(z_j-\bar z_j))^{m_j}$, which
violates (2.4).  Thus we must have $Q(0)=0$.\qed

\proclaim Proposition 2.5.   Let $Q\in\vf \mu$, $\mu\ne 0$, be
a vector field of the form (2.5).  If (2.6) holds, then $Q=0$.

\noindent{\sl Proof.}  Let us suppose that $\mu>0$.  (The case $\mu<0$ is
similar.)  The tangency condition is
$$\Re e\,Qp=\sum_{l=-L}^L\Re e\,Qp^{(\nu_l)}=0.\eqno(2.7)$$
We note that $Qp^{(\nu_L)}$ is the only term in (2.7) with signature $\mu+\nu_L$,
which is the largest possible, and there are no terms of signature $-\mu-\nu_L$
which might cancel with it upon taking the real part.  Thus
$$Qp^{(\nu_L)}=0.$$
The terms of signature $\nu_{L-1}+\mu$ in (2.7) must vanish, so this condition
is given either by 
$$Qp^{(\nu_{L-1})}$$ 
alone or by
$$Qp^{(\nu_{L-1})}+\overline{ Qp^{-(\nu_L)}}=0;\eqno(2.8)$$
the occurrence of the second case depends on whether
$\nu_{L-1}+\mu=-(-\nu_L+\mu)$ or not.

We show now that in  case (2.8)we have $Q^2p^{(\nu_{L-1})}=0$.  We recall that
$$p^{(-\nu_L)}=\overline{ p^{(\nu_L)} }=\sum \overline{f_{\nu_L,B}(z)}z^B$$
which gives
$$\overline {Qp^{(-\nu_L)}}=\sum f_{\nu_L,B}(z)\overline{Qz^B}.\eqno(2.9)$$
Since $Qp^{(\nu_L)}=0$, we have $Qf_{\nu_L,B}=0$ for all $B$.  Combining this
with (2.9), we have
$$Q\overline{Qp^{(-\nu_L)}} =\sum Qf_{\nu_L,B}\overline{Qz^B}=0.$$
From (2.8), then,
we have $Q^2p^{(\nu_{L-1})}=0$.  Continuing in this fashion, we have  $$Q^\alpha
p^{(\nu_l)}=0{\rm\ for\ some\ }\alpha\le L-l+1.$$

Now set $A:=\{Q=0\}$.  In case $A=\{0\}$, we let $S$ denote the orbit given by
Lemma 2.1.  By the remarks above, we must have $Q(0)=0$, so that by Lemma 2.3 it
follows from (2.4) that $p^{(\nu_l)}|_S=0$ for $-L\le l\le L$.  Thus we have
$S\subset\{p=0\}$ which contradicts (2.6).

The other case is $A\ne\{0\}$.  If $Q\ne0$, then we may let $\Sigma$ denote the
set given by Lemma 2.2.  We may assume that $\Sigma\subset\{Q=0\}$, for
otherwise if $c\in\Sigma\cap\{Q\ne0\}$, then as in the proof of Lemma 2.1, the
orbit $S$ passing through $c$ contains 0 in its closure, and we derive a
contradiction as in the previous case.  

Now we show that
$$p^{(\nu_l)}|_\Sigma=0 {\rm\ for\ all\ }-L\le l\le L.\eqno(2.10)$$
Comparing terms of signature $\nu_l$ in (2.7), we must have either
$Qp^{(\nu_l)}=0$ or
$$Qp^{(\nu_l)}+\overline{Q\overline{p^{(\nu_k)}}}=0\eqno(2.11)$$
if there is a $k$ such that  $\nu_l+\mu=-(-\nu_k+\mu)$.    From (2.9) and (2.11)
we obtain 
$$Qp^{(\nu_l)} =-\sum_B f_{\nu_k,B}(z)\overline{Qz^B}.\eqno(2.12)$$
Let $\kappa=(1-\nu_l)/2$ denote the weight of the indices $B$ appearing in (2.12),
and let ${\cal P}_\kappa(z)$ denote the holomorphic polynomials in $z$ of weight
$\kappa$.  Let $\varphi_1,\dots,\varphi_N$ denote a basis for the space ${\cal
P}_{\kappa+\mu}(z)/Q{\cal P}_\kappa(z)$.  It follows that
$\{z^I\bar\varphi_j\}$ forms a basis for the space
$${\cal P}_{\lambda,\kappa+\mu}(z,\bar z)/({\cal P}_\lambda(z)\overline{Q{\cal
P}_\kappa(z)})$$
where ${\cal P}_{\lambda,\kappa+\mu}(z,\bar z)={\cal
P}_\lambda(z)\overline{{\cal P}_{\kappa+\mu}(z)}$.  Since $\{\varphi_j\}$ is a
basis, there exist holomorphic polynomials $g_B$ and $g_j$ such that
$p^{(\nu_l)}=\sum g_B\overline{Qz^B}+\sum g_j\overline{\varphi_j}$.  This yields
$$Qp^{(\nu_l)}=\sum Qg_B\overline{Qz^B}+\sum
Qg_j\overline{\varphi_j}.\eqno(2.13)$$ The difference between (2.12) and (2.13) must
vanish, and $\{\varphi_j\}$ is a basis, so $\sum Qg_j\bar\varphi_j=0$.  Thus
$$Q\left(p^{(\nu_l)}-\sum g_B\overline{Qz^B}\right)=0.$$
By Lemma 2.2, we conclude that 
$$p^{(\nu_l)}-\sum g_B\overline{Qz^B}=0{\rm\ \ on\ }\Sigma.$$
Finally, $Q=0$ on $\Sigma$ since $\Sigma\subset A$, so 
$p^{(\nu_l)}|_\Sigma=0$.  Thus we  have $\Sigma\subset\{p=0\}$, which
contradicts (2.6).\qed

\proclaim Lemma 2.6.  If $\wt
Q\ge 0$, and $Q$ is nonzero, then $q_0$ is not identically zero and is divisible
by $w$.

\noindent{\sl Proof.}  If $q_0$ is not divisible by $w$, then it contains a
holomorphic polynomial in the variables $z_1,\dots z_n$.  Since $wt\,Q\ge0$, we
must have $wt\,q_j>0$ for all $j$.  Thus $q_j{\partial p\over\partial z_j}$ can
contain no pure holomorphic or antiholomorphic terms.  By equation (2.3), we
must have $\Re e\,iq_0=0$, so that $q_0=0$.  It follows, then, by Lemma 2.5 that
$Q=0$.
\qed

Next we classify the vector fields of negative weight.  It is
obvious that $\vf\mu$ is trivial if $\mu<-1$, and $\vf{-1}$ is generated by
${\partial\over\partial w}$.  It turns out that the only other possible weights
are the $-\delta_j$ themselves.  
\proclaim Lemma 2.7.  Let a nonzero vector field $Q\in\vf\mu$, $\mu<0$, be
given.  Then $\mu=-\delta_j$ for some $1\le j\le n$, and there is a weighted
homogeneous change of coordinates $(\tilde z_1,\dots,\tilde z_n)$, depending on
$(z_1,\dots,z_n)$ such that after possibly relabeling the indices $Q$ is given in
the new coordinates as $s_0{\partial\over\partial\tilde w} +
{\partial\over\partial\tilde z_1}$.  The new defining function $\tilde p$
satisfies (2.15) and (2.16) below, and 
$$s_0(z_1,0,\dots,0)=mc\left({z_1\over
2i}\right)^{m-1}$$ 
for some real constant $c$.

\noindent{\sl Proof.}  It is evident that if $wt\, Q<0$, then each coefficient
$q_j$ depends only on $z_1,\dots,z_n$, i.e. it is independent of $w=z_0$.  If
$q_0=0$, then $Q=0$ by Lemma 2.5.  Thus we suppose that $q_0\ne0$.  If $q_0$ is a
real constant, then $wt\,Q=-1$, and $Q=c{\partial \over \partial w}$, which
completes the proof in this case.  

Otherwise, $q_0$ is a holomorphic polynomial
of positive weight, which will need to be cancelled by pluriharmonic terms in
(2.3).  By (2.4) ${\partial p\over\partial z_j}$ cannot contain any purely
holomorphic monomials.  And if $q_j\ne0$ and $q_j(0)=0$ (i.e. if $wt\, q_j>0$),
then $q_j{\partial p\over\partial z_j}$ contains no pluriharmonic monomials,
so there is nothing pluriharmonic to cancel $q_0$ in (2.3).  Thus there must be
a $j$ such that $q_j(0)\ne0$, i.e. $q_j$ is constant, and thus $wt\,Q=-\delta_j$.

For convenience of notation, we may assume that $j=1$.  Let us define
homogeneous polynomials $h_2(z),\dots,h_n(z)$ by ${\partial h_j\over \partial
z_1}=q_j$.  It is immediate that $wt\,h_j=wt\,q_j+\delta_1$; and since
$wt\,Q=-\delta_1$, we have $wt \,q_j=-\delta_1+\delta_j$.  Thus
$wt\,h_j=\delta_j$, and the coordinate change defined by 
$$w=\tilde w,\quad z_1=\tilde z_1,\quad z_j=\tilde z_j+h_j(\tilde z),\quad 2\le
j\le n$$
is weighted homogeneous.  Since ${\partial\over\partial\tilde
z_1}=\sum{\partial z_j\over\partial\tilde z_1}{\partial\over\partial z_j}$, we
have $Q=s_0{\partial\over\partial\tilde w}+{\partial\over\partial\tilde z_1}$
if we define $s_0$ by $s_0(\tilde z)=q_0(z)$. 

We define a homogeneous polynomial $S$ by requiring that it be divisible by
$\tilde z_1$ and that ${\partial S\over\partial \tilde z_1}=s_0$.  We define
coordinates $(\hat w,\hat z)$ by setting
$$\tilde w=\hat w+S(\hat z),\quad \tilde z=\hat z.$$
The surface $\{\tilde p(\tilde z)+\tilde v=0\}$ becomes $\{\tilde p(\hat z)+\Im
m\,S(\hat z)+\hat v=0\}$.  In the $\hat{}$-coordinates we have
$Q={\partial\over\partial\hat z_1}$.  

Let us now drop the hats from the
coordinates.  The condition that $Q$ is tangential to
$\partial\Omega$ is equivalent to the condition that the function $\tilde
p(z)+\Im m\,S(z)$ is independent of the variable $\Re e\,z_1$.  Thus we may
write
$$\eqalign{ \tilde p(z)+{1\over 2i}(S(z)-\bar S(z))&=c\left({z_1-\bar z_1\over
2i}\right)^{m} +\cr
&+2\Re e\,\sum_{k=2}^n\alpha_k z_k\left({z_1-\bar z_1\over
2i}\right)^{m_k} +O\left(\sum_{k=2}^n|z_k|^2\right).\cr}$$ 
with $1\le m_k\le m$.  We claim that, in addition, we have 
$$m/2\le m_k\le
m-2.\eqno(2.14)$$
Since $p$ has weight 1, we must have $m\delta_1=1$ and
$\delta_k+m_k\delta_1=1$.  Since $\delta_k$ is the reciprocal of an integer, we
have $1-m_k/m\le{1\over 2}$, which gives the lower bound on $m_k$.  For the
upper bound, it is obvious that $m_k\le m-1$.  If equality holds, then
$\delta_1=\delta_k$, i.e. $z_1$ and $z_k$ have the same weight.  Thus the
coordinate change $\tilde z_1=z_1-i\alpha_k z_k$ is homogeneous, and in it the
monomial $z_ky_1^{m_k}$ does not appear.

By (2.4), the pure holomorphic terms are inside $S$, so
we have 
$$\tilde p(z_1,0,\dots,0)=c\left[\left({z_1-\bar z_1\over
2i}\right)^m-2\Re e\left({z_1\over 2i}\right)^m\right]\eqno(2.15)$$
$${\partial\tilde p\over\partial
z_j}(z_1,0,\dots,0)=\alpha_j\left[\left({z_1-\bar z_1\over
2i}\right)^{m_j}-\left({z_1\over 2i}\right)^{m_j}\right]\eqno(2.16)$$
and
$${1\over 2i}S(z_1,0,\dots,0)=c\left({z_1\over 2i}\right)^m.$$The derivative,
$s_0$, thus has the form stated above.
\qed

\bigskip
\centerline{\bf \S3.  Balanced Domains}
\bigskip

$\vf0$ contains the vector field
$${\cal D}=w{\partial\over \partial w} +\sum\delta_jz_j {\partial\over \partial
z_j},$$
which corresponds to homogeneous dilation.  For any vector field $Q\in\vf\mu$,
$[{\cal D},Q]=\mu Q$.  If $p$ is weighted homogeneous of weight $\mu$,
then the familiar Euler identity may be recast in the form $2\Re e\,{\cal
D}p=({\cal D}+\overline{{\cal D}})p=\mu p$. An analogue of this which will be
useful later is  $$2{\cal D}p=\mu p {\rm\ if\ 
and\ only\ if\ } p{\rm\ is\ balanced}.$$
If $p$ is not balanced, then we may consider $p^{(\nu)}$ as in (2.1).  In this
case we have
$${\cal D}p^{(\nu)}={\mu+\nu\over 2}p^{(\nu)}.\eqno(3.1)$$

\proclaim Lemma 3.1.  If $Q\in\vf0$, then $Q=c{\cal D}+{\cal L}$ for
some $c\in\RR$, and ${\cal L}$ is of the form (2.5).

\noindent{\sl Proof.}  By Lemma 4, we must have $q_0=cw$ for some constant
$c$.  It is evident that ${\cal L}:=Q-c{\cal D}\in \vf{0}$ and has the form
(2.5).
\qed

\proclaim Lemma 3.2.  If $Q\in\vf\mu$, $0<\mu<1$, and if $Q\ne0$, then
$\mu={1\over2}$.

\noindent{\sl Proof.}  Let us define $S:=[{\partial\over\partial
w},Q]=\sum{\partial q_j\over\partial w}{\partial\over\partial z_j}$.  Then
$S\in\vf{\mu-1}$, and by Lemma 2.6, $S\ne0$.  By Lemma 2.7, then,
$\mu-1=-\delta_j$ for some $j$.  In particular, if $\mu\ne{1\over2}$, then
${1\over2}<\mu=1-\delta_j<1$.  By Lemma 2.7, we may assume that
$S=s_0(z){\partial\over\partial w}+{\partial\over\partial z_1}$.  By Lemma 2.6,
$q_0$  is divisible by $w$, so
$$Q=ws_0{\partial\over\partial w} +(w+r_1(z)){\partial\over\partial z_1}
+\sum_{j=2}^nr_j(z){\partial\over\partial z_j}$$
for some homogeneous polynomials $r_j$.  

Now let us calculate the commutator
$$[S,Q]=\left (s_0^2-\sum r_j{\partial s_0\over\partial
z_j}\right){\partial\over\partial w}+\left(s_0+{\partial r_1\over\partial
z_1}\right){\partial\over\partial z_1} +\sum_{j=2}^n{\partial r_j\over\partial
z_1}{\partial\over\partial z_j}.$$
Since $\delta_1<{1\over2}$, the commutator has weight $1-2\delta_1>0$. 
Further, since the coefficient of ${\partial\over\partial w}$ is not divisible
by $w$ unless it is 0, we see that $[S,Q]=0$ by Lemma 2.6.  In particular,
$\partial r_j/\partial z_1=0$ for $2\le j\le n$.  If we set $z_2=\dots=z_n=0$,
then $r_j=0$ for $2\le j\le n$.  By homogeneity, $s_0(z_1,0,\dots,0)=\alpha
z_1^{m-1}$ and $r_1(z_1,0,\dots,0)=\beta z_1^m$. Since the coefficient of
${\partial \over\partial z_1}$ vanishes, we must have $\alpha+m\beta=0$.  But
since the coefficient of ${\partial\over\partial w}$ vanishes, too, we have
$\alpha^2-(m-1)\alpha\beta=0$.  Thus $\alpha=\beta=0$.  But by Lemma 2.7, we
have $\alpha=cm=0$, so $c=0$.  This contradiction proves the Lemma. \qed

We note that a polynomial $p$ can be nondegenerate in the sense of satisfying
(2.6), but $p^{(0)}$ may be degenerate.  We consider a different nondegeneracy
condition on a homogeneous polynomial $\varphi$:
$${\rm There\ is\ no\ holomorphic\ vector\ field\ }R\ne0{\rm\ such\ that\
}R\varphi=0.\eqno(3.2)$$
If $\varphi$ is strictly psh at some point, then (3.2) holds.

\proclaim Lemma 3.3.  Suppose that the balanced part $p^{(0)}$ of
$p$ satisfies (3.2), and let $\mu=\wt p$. If $l_k$, $1\le k\le n$ are weighted
homogeneous of weight $\delta_k$, and if $p+\sum_{k=1}^ml_k{\partial p\over
\partial z_k}=0$, then $l_k=-2\mu^{-1}\delta_k z_k$, and  $p$ is balanced.

\noindent{\sl Proof.}  By (3.1) the operator $L=\sum
l_k{\partial\over\partial z_k}$ preserves the splitting (2.1) in the sense
that $\delta(Lp^{(\nu)})=\nu$.  It follows that $p^{(0)}= -Lp^{(0)}$.  Since
$p^{(0)}$ is balanced, we have  
$$\sum_{k=1}^m(l_k+{2\mu^{-1}}\delta_kz_k){\partial
p^{(0)}\over \partial z_k}=0.$$  Thus the integral curves of the vector field
$L+{2\mu^{-1}}{\cal D}$ lie in the level sets of $p^{(0)}$, which contradicts
(2.6), unless $l_k=-2\mu^{-1}\delta_k z_k$.  It follows now from (3.1) that $p$
is balanced.
\qed

\proclaim Lemma 3.4.  Suppose that $p$ satisfies (2.6) and that $p^{(0)}$
satisfies (3.2).  If there is a nonzero vector field $Q\in\vf{{1\over 2}}$, then
$p$ is a balanced polynomial; and after a homogeneous change of coordinates and a
permutation of variables, $$Q=\lambda\left(-2iwz_1{\partial\over\partial w}
+w{\partial\over\partial z_1} -\sum_{j=1}^n2i\delta_j
z_1z_j{\partial\over\partial z_j}\right)\eqno(3.3)$$
for some $\lambda\in\RR$.

\noindent{\sl Proof.}  Let us use the notation $z_1,\dots,z_d$ for the
variables of weight ${1\over2}$, and let $\zeta_1,\dots,\zeta_e$ for the
variables with weight $<{1\over2}$.  Let $Q$ be a nonzero vector field of
weight ${1\over2}$, and let $Q^{(-{1\over2})}:=[{\partial\over\partial w},Q]$. 
By Lemma 2.7 we may assume that $Q^{(-{1\over2})}=-2iz_1{\partial \over\partial
w}+{\partial \over\partial z_1}$.  Thus by Lemma 2.5 we have 
$$Q=-2iz_1 w{\partial\over\partial w} +w{\partial\over\partial z_1}
+\sum_{j=1}^dq_j(z,\zeta){\partial\over\partial z_j} +\sum_{j=1}^e\tilde
q_j(z,\zeta){\partial\over\partial \zeta_j}.$$
We may write 
$$p=\sum z_j\bar z_j
+\sum (z_j\varphi_j(\zeta,\bar\zeta)+\bar z_j\bar
\varphi_j(\zeta,\bar\zeta))+\tilde p(\zeta,\bar\zeta).\eqno(3.4)$$
After a change of coordinates of the form $z_j\mapsto z_j+\psi_j(\zeta)$, we
may assume that $\varphi_j$ contains no anti-holomorphic terms.  Then after
$w\mapsto w+\chi(\zeta)$, we may assume that $\varphi_j$ contains no
holomorphic terms, i.e.~that (2.4) is satisfied.    Thus (2.3) takes the form
$$\eqalign{ u\,\Re e\varphi_1 +\Re e&\left[-2iz_1v+iv\varphi_1 +\sum_{j=1}^d
q_j(\bar z_j+\varphi_j)\right.+\cr 
& \left.+\sum_{k=1}^e\tilde q_k\left( \sum_{j=1}^d
z_j{\partial\varphi_j\over\partial\zeta_k} +\sum_{j=1}^d \bar
z_j{\partial\bar\varphi_j\over\partial\zeta_k} +{\partial\tilde
p\over\partial\zeta_k}\right)\right]=0.\cr}\eqno(3.5)$$
We note immediately that the coefficient of $u$ must vanish, i.e.~that $\Re
e\,\varphi_1=0$.

Now
we set $u=0$ and $\zeta=0$ in equation (2.3).  Since $\varphi_j$ has weight 1,
${\partial\varphi_j(0)\over\partial\zeta_k}=0$, so we have  
$$\Re e(2iz_1\sum
z_j\bar z_j+\sum q_j(z,0)\bar z_j)=\Re e\sum \bar z_j(2iz_1z_j+q_j(z,0))=0.$$
We conclude, then, that 
$$q_j(z,\zeta)=-2iz_1z_j+ \sum z_k q^{(k)}_j(\zeta) +q_j^\prime(\zeta)$$
for $1\le j\le d$.  Similarly, since $\tilde q_j$ has weight $<1$, we may write
$$\tilde q_j=\sum z_k\tilde q^{(k)}_j(\zeta) +\tilde q^\prime_j(\zeta).$$

Now we observe that there are no coefficients (i.e.~functions of $\zeta$ and
$\bar\zeta$) of $\bar z_1^2$ inside the term in square brackets in (3.5). 
Thus the coefficient of $z_1^2$ must vanish, i.e.
$$\sum_{k=1}^e \tilde q^{(1)}_k{\partial\varphi_1\over\partial\zeta_k}=0.$$
Similarly, taking the coefficient of $z_1\bar z_1$ we have
$$\Re e\left(q_1^{(1)}+ i\varphi_1+\sum_{k=1}^e \tilde
q^{(1)}_k{\partial\bar\varphi_1\over\partial\zeta_k}\right)=0.$$
The only purely holomorphic or antiholomorphic terms come from $q^{(1)}_1$ so
$q^{(1)}_1=0$.  Now since $\varphi_1$ is pure imaginary, these two equations give
$\varphi_1=0$. 

There is no function of $\zeta,\bar\zeta$  as multiple of $\bar z_1\bar z_m$ for
$2\le m\le d$ in the bracketed term in (3.5).  Thus the coefficient of $z_1z_m$
must vanish:
$$\sum_{k=1}^e\tilde q^{(1)}_k{\partial\varphi_m\over\partial\zeta_k}=0.$$
The coefficient of $z_1\bar z_m$ plus the conjugate of the coefficient of
$z_m\bar z_1$ must also vanish, so
$$q_m^{(1)} +\bar q_1^{(m)} +\sum_{k=1}^e\tilde
q^{(1)}_k{\partial\bar\varphi_m\over\partial\zeta_k}=0.$$
It follows upon comparing pure terms that $q_m^{(1)}=q_1^{(m)}=0$.

The coefficient of $\bar z_1$ in the bracketed term in (3.5) is
$q^\prime_1$.  Adding the conjugate of this to the coefficient of $z_1$, we
obtain $$\bar q^\prime_1 + 2i\tilde p +\sum_{k=1}^e \tilde
q^{(1)}_k{\partial\tilde p\over\partial\zeta_k}=0.$$
The only pure terms come from $\bar q_1^\prime$, so we must have
$q_1^\prime=0$.  By Lemma 3.3, then, we conclude that
$$\tilde q^{(1)}_k=-4i\delta_k\zeta_k\eqno(3.6)$$
and that $\tilde p$ is balanced.

We note that there are no terms of the form $\bar z_1\bar z_j$ in the bracketed
term in (3.5) for $2\le j\le d$, so the coefficient of $z_1z_j$ must vanish,
i.e. 
$$\sum_{k=1}^e q^{(1)}_k {\partial\varphi_j\over\partial\zeta_k}=0.$$
But by (3.6), this gives $\varphi_j=0$, since $\varphi_j$ contains no purely
anti-holomorphic terms.

Now we may inspect the coefficients of $z_j\bar z_m$ for $2\le j,m\le d$, and
since $\varphi_j=0$, we get $q^{(m)}_j=0$.

Finally, the coefficients of $z_j$ and $\bar z_j$ for $2\le j\le d$ give
$$\bar q^\prime_j +\sum_{k=1}^e\tilde q^{(j)}_k{\partial\tilde
p\over\partial\zeta_k}=0.$$
Again, the pure terms vanish, so $q^\prime_j=0$.  By Proposition 2.5, we have
$\tilde q^{(j)}_k=0$.
Setting $z=0$ in (3.5), we have
$$\Re e\sum_{k=1}^e\tilde q^\prime_k{\partial\tilde p\over\partial\zeta_k}=0,$$
so that $\tilde q^\prime_k=0$ by Proposition 2.5.  This completes the proof.
\qed

\noindent{\bf Remark.}  We have in fact shown that under the hypotheses of
Lemma 3.4 we have
$$p=\sum_{j=1}^d z_j\bar z_j +\tilde p(\zeta,\bar\zeta).$$

\proclaim Lemma 3.5. Suppose that $p$ satisfies (2.6), and that $p^{(0)}$, and
$(p^2)^{(0)}$ satisfy (3.2).  If there is a nonzero vector field $Q\in\vf{1}$,
then $p$ is a balanced polynomial; and after a homogeneous change of
coordinates and a permuation of variables, 
$$Q=\lambda\left(w^2{\partial\over\partial
w}+\sum_{j=1}^n2\delta_jwz_j{\partial\over\partial z_j}\right)\eqno(3.7)$$ 
for
some $\lambda\in\RR$.

\noindent{\sl Proof.}  By Lemmas 2.6 and  2.7, $q_0$ is a real multiple of
$w^2$.  Thus
$$Q=w^2{\partial\over\partial w} +
\sum_{j=1}^n(wq_j(z)+r_j(z)){\partial\over\partial z_j}$$
where $wt\,q_j=\delta_j$, and $wt\,r_j=1+\delta_j$.  The coefficient of $u$ in
(2.3) is then
$$\Re e\left(-p+\sum_{j=1}^nq_j{\partial p\over\partial z_j}\right)=0.$$
Multiplying by $2p$, we have
$$\Re e\left(-2p^2+\sum_{j=1}^nq_j{\partial (p^2)\over\partial
z_j}\right)=0.\eqno(3.8)$$

Now we set $u=0$ in (2.3) and obtain
$$\Re e\left[-i\sum pq_j{\partial p\over\partial z_j} +\sum r_j{\partial
p\over\partial z_j}\right]=0$$
or, after doubling,
$$\Im m\left(\sum q_j{\partial(p^2)\over\partial z_j}\right) +2\,\Re
e\left(\sum r_j{\partial p\over\partial z_j}\right)=0.\eqno(3.9)$$
Every monomial in $r_j{\partial p\over\partial z_j}$ is of the form $z^A\bar
z^B$ with $wt\,A\ge 1+\delta_j$ and $wt\,B\le 1-\delta_j$.  Thus the second
term in (3.9) can have no balanced monomials.  The operator $\sum
q_j{\partial\over\partial z_j}$ preserves balanced polynomials, so we may add
(3.8) and (3.9) and take the balanced part to obtain
$$\sum_{j=1}^nq_j{\partial(p^2)^{(0)}\over\partial z_j}=2(p^2)^{(0)}$$
where $(p^2)^{(0)}$ denotes the balanced part of $p^2$.  By Lemma 3.3, we have
$q_j=2\delta_jz_j$.

Now we add (3.8) and (3.9) to obtain
$$2\cD(p^2)=\sum_{j=1}^n2\delta_jz_j{\partial(p^2)\over\partial z_j}=2p^2-2i\,\Re
e\,\sum_{j=1}^nr_j{\partial p\over\partial z_j}.\eqno(3.10)$$
Now we write 
$$p^2=(p^2)^{(-\mu_1)}+\dots+(p^2)^{(\mu_1)}$$
as in (2.1).  By (3.1) we have
$$\sum_{j=1}^n2\delta_jz_j{\partial(p^2)^{(\mu_j)}\over\partial
z_j}=\sum_{j=1}^n(2+{\mu_j})(p^2)^{(\mu_j)}.$$
Thus if we can show that $r_j=0$ for $1\le j\le n$, then from (3.10) we will
have 
$$(1+{\mu_j\over 2})(p^2)^{(\mu_j)}=(1-{\mu_j\over
2})\overline{(p^2)^{(-\mu_j)}}=(1-{\mu_j\over
2})(p^2)^{(\mu_j)}.$$
We conclude that the only terms that can appear correspond to $\mu_j=0$.  Thus
$p^2$ is a balanced polynomial.  It follows that $p$ is balanced, too. With
Lemma 3.6, then, the proof will be complete. \qed

\proclaim Lemma 3.6.  With the notation of Lemma 3.5, $R:=\sum
r_j{\partial\over\partial z_j}=0$.

\noindent{\sl Proof}.  From (3.1) and (3.10) we obtain
$i\mu(p^2)^{(\mu)}=R p^{(\mu-1)}$
for $\mu>0$, or
$$i(\nu+1)(p^2)^{(\nu+1)}=Rp^{(\nu)}\eqno(3.11)$$
for all indices $\nu=\nu_{-L},\nu_{-L+1},\dots,\nu_L$.

The largest signature in $(p^2)^{(\nu)}$
is $2\nu_L<1+\nu_L$, so setting $\nu=\nu_L$ in (3.11) we obtain
$$Rp^{(\nu_L)}=0.$$
Now proceed by induction to show that
$$R^{2^j}p^{(\nu_{L-j})}=0{\rm\ for\ }0\le j\le 2L.\eqno(3.12)$$
We consider (3.11) 
for $\nu=\nu_j$.  Each element of $(\nu+1)(p^2)^{(\nu+1)}$ is made up from a
product of the form $(\nu+1)p^{(\nu_a)}p^{(\nu_b)}$ with $\nu_a+\nu_b=\nu_j+1$. 
Since $\nu_a$, $\nu_b<1$, it follows that $\nu_a,\nu_b>\nu_j$. 
By induction, $R^{2^{j-1}}p^{(\nu)}=0$ for $\nu=\nu_a$ and $\nu=\nu_b$; and so 
$$R^{2^j-1}(p^{(\nu_a)}p^{(\nu_b)})=0.$$
This proves (3.12).

Next we consider the ordered family ${\cal
H}=\{h_1,h_2,\dots,h_N\}$ of homogeneous polynomials defined as follows.  We
list out ${\cal H}$ in groups in the order ${\cal G}_L,{\cal
G}_{L-1},\dots,{\cal G}_{-L}$.  For fixed $l$, $-L\le l\le L$ we let
$B_1,B_2,\dots$ be an arbitrary ordering of the (finitely many) multi-indices
appearing in $p^{(\nu_l)}=\sum f_{\nu_l,B}\bar z^B$.  Then we define ${\cal
G}_l$ to be the (finite) ordered set
$$\{R^{2^{L-l}-1}f_{\nu_l,B_1},R^{2^{L-1}-2}f_{\nu_l,B_1},\dots,f_{\nu_l,B_1},
R^{2^{L-l}-1}f_{\nu_l,B_2},R^{2^{L-l}-2}f_{\nu_l,B_2},
\dots,f_{\nu_l,B_2},\dots\}.$$
Let us define $V_0=\cx n$ and $V_m=\{h_1=\dots=h_m=0\}$ for $1\le m\le N$. 
Clearly $V_0\supset V_1\supset\dots\supset V_N$.  Since each $f_{\nu_l,B}$
vanishes on $V_N$, it follows that $V_N\subset\{p=0\}$.  Thus dim$V_N=0$. 
And since $V_N$ is invariant under $D_\tau$ for all $\tau\in\cx{}$, we
have $V_N=\{0\}$. 
Further, by the choice of ordering of ${\cal
H}$, thogether with (3.12), we see that 
$$Rh_{j+1}=0{\rm\  on\ }V_j{\rm\ for\ every\ }j.\eqno(3.13)$$  
Thus
$R$ is tangent to $V_m$ at the regular points of $V_m$ where $R\ne0$. 

Since we are trying to prove that $R=0$, we may assume that $A:=\{R=0\}\ne\cx
n$.  Thus we may choose $m$ such that $V_{m+1}\cap A$ has a component which has
codimension 1 in $V_m$.  Passing to irreducible components, we may assume that
$V_m$ and $V_{m+1}$ are irreducible, and $V_{m+1}\subset A$.  There are now two
cases to consider.  The first case is that $V_{m+1}=\{0\}$.  Thus $R\ne0$ on
$V_m-\{0\}$, and by Lemma 2.1 there is an orbit $S$ of $R$ with 0 in its
closure.  It follows from (3.12) and Lemma 2.3 that $p^{(\nu)}=0$ on $S$.  Thus
$p=0$ on $S$, which contradicts (2.6).

In the other case, $dim(V_m)\ge 1$, and we may choose a (constant) tangent
vector $T=\sum\alpha_j{\partial\over\partial z_j}$ which is tangent to $V_m$ at
some point, and we define $\varphi:=T^kh_{m+1}$.    We may choose $k$ such
that $V_{m+1}\subset\{\varphi=0\}\cap V_m$ and $d\varphi|_{V_m}$ does not
vanish identically on $V_{m+1}$.  Then there is a holomorphic vector field
$\tilde R$ on $V_m$ such that $\varphi^d\tilde R=R$ for some $d\ge 1$ and such
that $\tilde R$ does not vanish identically on $V_{m+1}$.  By (3.13) we have
$\varphi^d\tilde R h_{m+1}=Rh_{m+1}=0$ on $V_m-V_{m+1}$.  Thus $\tilde
Rh_{m+1}=0$ on $V_{m+1}$.  It follows that $\tilde R$ is tangential to the
regular points of $V_{m+1}$ on the (nonempty) set where $\tilde R\ne0$.

For each $\mu$, we let $s_\mu\ge0$ denote the largest integer such that 
$$p^{(\mu)}=O(|\varphi|^{s_\mu})$$
holds on $V_m$.  By (3.11) and the definition of $\tilde R$ it follows that
$s_\mu\ge d$ for $\mu>0$.  Since $p^{(-\mu)}=\overline{p^{(\mu)}}$, we have
$s_\mu\ge d$ for $\mu<0$.  Let $s$ be the minimum value of $s_\mu$ for
$\mu\ne0$.  Since $(p^2)^{(1)}$ consists of products $p^{(\nu_a)}p^{(\nu_b)}$
with $0<\nu_a,\nu_b<1$ and $\nu_a+\nu_b=1$, it follows that 
$$(p^2)^{(1)}=O(|\varphi|^{2s}).$$

From (3.10) we have $Rp^{(0)}=i(p^2)^{(1)}$, so we have
$$\varphi^d\tilde Rp^{(0)}=O(|\varphi|^{2d}).$$
It follows, then, that $\tilde Rp^{(0)}=0$ on $V_m$.  By Lemma 2.2, then,
$p^{(0)}$ vanishes on a variety passing through 0.  Thus $p$ vanishes on the
same variety, which contradicts (2.6).\qed

We may summarize the work of \S3 by the following.

\proclaim Theorem 3.7.  Let $p$ be homogeneous of weight 1, let $p$ satisfy
(2.6), and let $p^{(0)}$ and $(p^2)^{(0)}$ both satisfy (3.2).  If there is a
tangential holomorphic vector field $Q$ for $\{v+p(z)<0\}$ with 
$\wt Q>0$, then $p$ is balanced.

\bigbreak
\centerline{\bf\S4.  Domains with Noncompact Automorphism Groups}
\bigskip

In this Section we give the proofs of Theorems 1 and 2.
\bigskip

\noindent{\sl Proof of Theorem 1.}  Let $p$ be homogeneous.  We may obtain the
balanced part $p^{(0)}$ as
$$p^{(0)}(z) = {1\over2\pi}\int_0^{2\pi}
p(e^{i\theta}z_1,\dots,e^{i\theta}z_n)d\theta.$$
If $p\ge0$, then $p^{(0)}\ge0$ and $(p^2)^{(0)}\ge0$.  Since $p^{(0)}$ is
symmetrized, it follows that $p^{(0)}$ is invariant under ${\cal D}_\tau$ for
all $\tau\in\cx{}$.  Thus, for $z_0\ne0$, $p^{(0)}(z_0)>0$, since otherwise
$p^{(0)}$ (and thus $p$) would vanish on the ${\cal D}_\tau$-orbit of $z_0$. 
It follows that the level sets $\{p^{(0)}=c\}$ are compact, and thus there are
points where $p^{(0)}$ is strongly psh.  So (3.2) holds.  Similarly, (3.2)
holds for $(p^2)^{(0)}$.

Now let $Q$ denote the homogeneous part of $H$ of lowest weight.  If $Q$ has
weight $\mu>1$, then the commutator $[Q,{\partial\over\partial w}]$ has weight
$\mu-1$ and is not the zero vector field by Lemma 2.6.  Taking further
commutators, we may assume that $Q$ has weight $0<\mu\le 1$.  By Lemma 3.2, then
$\mu$ is either ${1\over2}$ or 1.  Applying Lemmas 3.4 and 3.5, we may assume
that $Q$ has the form of either (3.3) or (3.7).  In either case, we may apply the
argument of the final rescaling in \S5 of [BP] to conclude that $\Omega$ is
biholomorphically equivalent to $\{v+p<0\}$.\qed

\noindent{\sl Admissible Assignment of Weights}.
For the rest of this Section we let $\Omega$ denote a smooth, convex surface with
finite type boundary.  Let us fix a point $0\in\partial\Omega$, and assume that
the tangent plane to $\partial\Omega$ at 0 is given by $\{v=0\}$.  Let us begin
by showing how to make an admissible assignment of weights.  We assign weight 1
to the variable $w$, and we write $\partial\Omega$ as $\{v+f(u,z)=0\}$ in a
neighborhood of 0. For a tangent vector $T=\sum a_j{\partial\over\partial
z_j}$, we let $Ord(f(0,z),T)=m$ be the smallest positive integer such that
$T^k(f(0,0))=0$ for $1\le k\le m-1$ and $T^mf(0,0)\ne0$.  Since $\Omega$ has
finite type, each vector $T\ne0$ has a finite order.  We let $L_1$ be the set
of complex tangent vectors such that $Ord(f(0,z),T)=m_1$ is maximum.  Then we
define numbers $m_1>\cdots>m_k$ and complex subspaces $L_1\subset L_2\subset
\dots\subset L_k=\cx n$ by the condition that $Ord(f(0,z),T)=m_j$ for $T\in
L_j-L_{j-1}$.  After a complex linear change of coordinates, we may assume that
there are integers $n_1<n_2<\cdots<n_k$ such that the coordinate system
$\{z_1,\dots,z_n\}$ has the property that $L_j$ is spanned by
$\{{\partial\over\partial z_s}:n_{j-1}+1\le s\le n_j\}$.  We assign weight
$m_j^{-1}$ to the variables $\{z_s:n_{j-1}+1\le s\le n_j\}$.  By the convexity
of $f$, the $m_j$ are all even.

Now let $p$ denote the terms of weight 1 in the Taylor expansion of $f$ at
$z=0$.  We must show that all monomials in the Taylor expansion of $f-p$ at the
origin have weight greater than or equal to one.  If not, there is a term of
minimal weight $\mu<1$ in $f$.  Let $q(z)$ denote the terms of weight $\mu$.  If
we perform the scaling of coordinates $\chi_t(w,z_1,\dots,z_n)=(t^\mu
w,t^{\delta_1}z_1,\dots,t^{\delta_n}z_n)$, then the surface $\Omega$ is
transformed to the surface $\Omega_t$, which converges to $\{v+q(z)<0\}$ as
$t\to0$.  Since $\Omega$ is convex, it follows that $q$ is a convex function. 
By the construction of the weights, however, there are no monomials of the form
$x_j^{m}$ or $y_j^m$ appearing in $q$.  By convexity, then, $q=0$.  We conclude
that this assignment of weights is admissible.

\medskip

The polynomial $p$ itself is obtained as the limit of $f$ under the scalings
above with $\mu=1$.  Thus $p$ is convex, and $\{p=0\}$ is a real linear subspace
of $\cx n$.  Since the order of $p$ in any direction $T$ is finite, $p$ cannot
vanish on a complex line.  Thus $\{p=0\}$ is a totally real linear subspace of
$\cx n$, and hence (2.6) is satisfied.

We will invoke several results which were
proved under slightly different hypotheses in [BP].  In fact, since $\Omega$ is
convex, several technicalities in [BP] can be avoided here.
First, we apply Lemma 2 of [BP]: if $Aut(\Omega)$ is noncompact, then $\Omega$ is
biholomorphically equivalent to a domain $D=\{v+\rho(z_1,\dots,z_n)<0\}$, where
$\rho$ is a convex polynomial (not necessarily homogeneous).  By
construction, $\rho$ satisfies (2.4).  Let $g:D\to\Omega$ denote this
biholomorphism, and let
$H=g_*(2{\partial\over\partial w})$.  Since $\Omega$ is convex and finite type,
$H$ extends smoothly to $\bar\Omega$ and induces a tangential holomorphic vector
field.  The mapping $g$ extends to a homeomorphism between $\bar
D\cup\{\infty\}$ and $\bar\Omega$.  We translate coordinates so that
$g(\infty)=0\in\partial\Omega$.

Since $\Omega$ is convex, we assign weights as above, and let $Q$ denote the
part of $H$ with minimal weight.  Since $H$ vanishes to finite order, it
follows that $Q\ne0$.  Since $Q(0)=0$, it follows from Lemma 2.7 that $Q$ has
weight $\ge0$.  If $wt\,Q=0$, then by Lemma 3.1, $Q=c{\cal D} + {\cal L}$. 
However, if $c\ne0$, then 0 is either a  source or a sink, and in neither case
can it be parabolic.

\noindent{\sl Uniform Hyperbolicity}.  We will consider a family
of convex domains in $\cx {n+1}$, and we wish to have an estimate for the
Kobayashi metric $F$, which is uniform over the family of domains. We let 
$F(q,\xi,\Omega)$ denote either the Kobayashi metric for the domain 
$\Omega$ at a point $q\in\Omega$, in the direction $\xi\in T_q(\Omega)$.  

We consider two sorts of convergent domains.  The first is a
sequence of domains of the form $G^{(\nu)}=\{v+\rho^{(\nu)}(z)<0\}$, where
$\rho^{(\nu)}$ is a sequence of convex polynomials of degree $\le d$, which
converge to a convex polynomial $\rho^{(\infty)}$ which is nondegenerate in the
sense that $\{z:\rho^{(\infty)}(z)=0\}$ is a totally real linear subspace of
$\cx{n+1}$.  It follows that there is a nondegenerate convex polynomial
$\hat\rho$ such that $\rho^{(\nu)}\ge\hat\rho$, and so we can use $\hat
G=\{v+\hat\rho<0\}$ to obtain the  estimate $F(q,\xi,G^{(\nu)})\ge
F(q,\xi,\hat G)$  for all $\xi$ and $q\in G^{(\nu)}\cap \hat G$.  Since the
limit domain $\hat G$ is Kobayashi hyperbolic, it follows that for a compact
$K\subset\hat G$, there exists $\epsilon>0$ such that
$$F(q,\xi,G^{(\nu)})\ge\epsilon|\xi|\eqno(4.1)$$ holds for all $q\in K$ if $\nu$
is sufficiently large.

The second is a family of scalings $\Omega^{(\nu)}:=\chi^{(\nu)}(\Omega)$.  As was observed
above, the domains $\Omega^{(\nu)}$ converge to the domain $\Omega_{hom}$.  Further,
there exists an $\epsilon>0$ such that for any compact $K\subset\Omega_{hom}$
$$F(q,\xi,\Omega^{(\nu)})\ge\epsilon|\xi|\eqno(4.2)$$ holds for $q\in K$ and $t$
sufficiently large.

\noindent{\sl Erratum}.  Similar estimates were used in the proof of Lemma 7 in
[BP, p. 181].  It was incorrectly stated there that $\tilde\rho^*$ was strictly
psh on $\{z_1\ne0\}$.  However, we may set $c=0$, so that $\tilde\rho^*$ is
strictly psh on $\{\psi_{z_1}\ne0\}$, and the proof of Lemma 7 works without
change if the cases $\gamma=0$ and $\gamma\ne0$ are considered separately.  The
Remark after Lemma 7 is incorrect.  The domain $\{v+P(z_1,\bar
z_1)+z_2\bar z_2+\dots+z_n\bar z_n<0\}$ is in fact Carathedory complete (which
follows from [BF]).  But the weaker property of being Kobayashi hyperbolic is
not as trivially proved as was asserted in [BP].

\proclaim Lemma 4.1.  There is a biholomorphic mapping $h:D\to\Omega_{hom}$
such that $h_*({\partial\over\partial w})=cQ$ for some real number $c\ne0$. 
Further, $Q$ has weight $>0$.

\noindent{\sl Proof.}  We let $g:D\to\Omega$ be as above, and we consider a
sequence of mappings $h^{(\nu)}:G^{(\nu)}\to\Omega_\nu$ defined as follows.  We
let $\Omega_\nu:=\chi_\nu\Omega$, and we let
$q^{(\nu)}=(w^{(\nu)},z_1^{(\nu)},\dots,z_n^{(\nu)}):=g^{-1}\chi_\nu^{-1}(w^0,z^0)$,
for some point $(w^0,z^0)\in\Omega$.  Now we make a coordinate change 
$$\hat z_j=z_j-z_j^{(\nu)},\qquad1\le j\le n$$
$$\hat w=w-\Re e\,w^{(\nu)}-a_0^{(\nu)}i+\Re
e\sum a_j^{(\nu)}(z_j-z_j^{(\nu)}).$$
We follow this with a scaling of coordinates 
$$\hat
w=\lambda^{(\nu)}\tilde w,\quad \hat z_j=\mu_j^{(\nu)}\tilde z_j,\quad 1\le
j\le n.$$
Thus the domain $D$ takes the form $G^{(\nu)}=\{v+\rho^{(\nu)}(z)<0\}$.
We choose the $\lambda^{(\nu)}$ and $\mu^{(\nu)}$ such that the coordinates of
$q^{(\nu)}$ are $(-i,0,\dots,0)$, and for $1\le j\le n$ the modulus of
the largest coefficient of $\rho^{(\nu)}(0,\dots,z_j,\dots,0)$ is 1.  The
mapping $h^{(\nu)}$ is then defined by applying the change of coordinates which
takes $G^{(\nu)}$ to $D$ and  following this with $\chi_\nu\circ g$.

The $G^{(\nu)}$ are a family of convex domains which, if we pass to a
subsequence, will converge to a domain $G:=\{v+\rho^{(\infty)}(z)<0\}$.  By the
uniform hyperbolicity condition (4.2), $\{h^{(\nu)}\}$ is a normal family,
with a limit function $h:G\to\overline{\Omega}_{hom}$.  However, since
$h^{(\nu)}(-i,0,\dots,0)=(w^0,z^0)\in\Omega_{hom}$ and since we have estimate
(4.1), it follows that $h$ is a biholomorphism.

Finally, we observe that the change of coordinates on the domain $D$ dilates 
the vector field ${\partial\over\partial w}$ by a factor of
$1/\lambda^{(\nu)}$.  Applying $g_*$, we obtain a scalar multiple of the
vector field $H$ on $\Omega$.  Finally, the scalings ${\chi_\nu}_*$, applied to
a multiple of $H$, converge to a multiple of $Q$ as $\nu\to\infty$.  Since $h$
is a biholomorphism, we must have $c\ne0$.

Now we wish to show that $wt\,Q>0$.  By the remarks above, $wt\,Q>0$; and if
$wt\,Q=0$, the $Q={\cal L}$, as in Lemma 3.1.  But in this case, the orbits
of $\Re e{\cal L}$ lie inside the level sets $\{p=c\}$.  Since such
orbits do not occur for a parabolic fixed point, we must have $wt\,Q>0$. \qed

\bigskip
\noindent{\sl Proof of Theorem 2.}  Applying Lemma 4.1, we have
$h:G\to\Omega_{hom}$.  Since ${\partial\over\partial w}$ is a parabolic
vector field on $G$, it follows that $Q$ must be a parabolic vector field on
$\Omega_{hom}$ with weight $>0$.  By Theorem 1,
then, it follows that $\Omega$ is biholomorphically equivalent to a domain of
the form (1.2).\qed

\bigskip
\bigbreak
\centerline{\bf References}
\medskip

\item{[BF]}  E. Bedford and J.-E. Forn\ae ss,  A construction of peak functions
on weakly pseudoconvex domains.  Ann. Math. 107 (1978), 555--568.
\item{[BP]}  E.
Bedford and S. Pinchuk, Domains in ${\bf C}^{n+1}$ with noncompact automorphism
group, Jour. of Geometric Analysis, 1 (1991), 165--191.

\item{[F]}  S. Frankel, Complex geometry of convex domains that cover varieties, Acta. Math. 163
(1989), 109--149.

\item{[KU]} W. Kaup and H. Upmeier, On the automorphisms of circular and
Reinhardt domains in complex Banach spaces,  Manuscripta Math. 25 (1978),
97--133.

\item{[Ki]} K.-T. Kim, Complete localization of domains with noncompact
automorphism groups, Trans. A.M.S. 319 (1990), 130--153.

\bigskip

\noindent Indiana University, Bloomington, IN 47405
\bigskip

\noindent Bashkirian State University, Ufa, USSR
\bye